\def\updated{20 February 2006}

\count100= 53 \count101= 8

\message{: version 12dec05}

\font\titlefont=cmr17
\font\twelverm=cmr12
\font\ninerm=cmr9
\font\sevenrm=cmr7
\font\sixrm=cmr6
\font\fiverm=cmr5
\font\ninei=cmmi9
\font\seveni=cmmi7
\font\ninesy=cmsy9
\font\sevensy=cmsy9

\font\sixbf=cmbx6
\font\fivebf=cmbx5
\font\ninebf=cmbx9
\font\nineit=cmti9
\font\ninesl=cmsl9
\font\nineex=cmex9

\def\ninepoint{\def\rm{\fam0\ninerm}
    \textfont0 = \ninerm
    \textfont1 = \ninei
    \textfont2 = \ninesy
    \textfont3 = \nineex
    \scriptfont0 = \sevenrm
    \scriptfont1 = \seveni
    \scriptfont2 = \sevensy
    \scriptscriptfont0 = \fiverm
    \scriptscriptfont1 = \fivei
    \scriptscriptfont2 = \fivesy
    \textfont\itfam=\nineit \def\it{\fam\itfam\nineit}
    \textfont\bffam=\ninebf \scriptfont\bffam=\sixbf
    \scriptscriptfont\bffam=\fivebf \def\bf{\fam\bffam\ninebf}
    \textfont\slfam=\ninesl \def\sl{\fam\slfam\ninesl}
    \baselineskip 10pt}

\magnification\magstephalf

\hsize6.5truein\vsize8.6truein
\voffset.5truein


\def\title#1{\toneormore#1||||:}
\def\titexp#1#2{\hbox{{\titlefont #1} \kern-.25em%
  \raise .90ex \hbox{\twelverm #2}}\/}
\def\titsub#1#2{\hbox{{\titlefont #1} \kern-.25em%
  \lower .60ex \hbox{\twelverm #2}}\/}

\def\author#1{\bigskip\bigskip\aoneormore#1||||:\smallskip\centerline{\updated}}

\def\abstract#1{\bigskip\bigskip\medskip%
 {\ninepoint
 \narrower{\bf Abstract.~}\rm#1\bigskip}\bigskip}

\def\toneormore#1||#2||#3:{\centerline{\titlefont #1}%
    \def\next{#2}\ifx\next\empty\else\medskip\toneormore#2||#3:\fi}
\def\aoneormore#1||#2||#3:{\centerline{\twelverm #1}%
    \def\next{#2}\ifx\next\empty\else\smallskip\aoneormore#2||#3:\fi}

\def\tocline#1#2#3{\nexttoc{\noexpand\noexpand\noexpand\\\hskip#2truecm #1.  #3.  }}

\def\footnoterule{\kern -3pt \hrule width 0truein \kern 2.6pt}
\def\leftheadline{\ifnum\pageno=\count100 \hfill%
  \else\hfil\it\shortauthor\hfil\llap{\rm\folio}\fi}
\def\rightheadline{\ifnum\pageno=\count100 \hfill%
  \else\hfil\it\shorttitle\hfil\llap{\rm\folio}\fi}

\def\References{\goodbreak\bigskip\centerline{\bf References}%
   \bigskip\frenchspacing}

\nopagenumbers
\headline{\ifodd\pageno\rightheadline \else\leftheadline\fi}
\footline{\hfil}
\null\vskip 18pt
\centerline{}
\pageno=\count100
\count102=\count100
\advance\count102 by -1
\advance\count102 by \count101


\def\copyright{\hbox{{\twelverm o}\kern-.61em\raise .46ex\hbox{\fiverm c}}}

\insert\footins{\sixrm
\medskip
\baselineskip 8pt
\leftline{Surveys in Approximation Theory
  \hfill {\ninerm \the\pageno}}
\leftline{Volume 2, 2006.
pp.~\the\pageno--\the\count102.}
\leftline{Copyright \copyright\ 2006 Surveys in Approximation Theory.}
\leftline{ISSN 1555-578X}
\leftline{All rights of reproduction in any form reserved.}
\smallskip
\par\allowbreak}


\newcount\blackmarks\blackmarks0
\newcount\eqnum
\newcount\labelnum

\def\singlecount{\let\labelnum\eqnum}



\title{Discrete Linear Interpolatory Operators}
\author{J.~Szabados}

\def\shorttitle{Discrete Linear Interpolatory Operators}
\def\shortauthor{J.~Szabados}

\abstract{This is a survey
on discrete linear operators which, besides approximating in
Jackson or near-best order, possess some interpolatory property at
some nodes. Such operators can be useful in numerical
analysis.}

A central problem in approximation theory
is the construction of simple functions that approximate well a
given set of functions. Traditionally, by ``simple" we mean
polynomials or rational functions as they are easily implemented
on computers, and the set of functions is generally characterized
by continuity or belonging to some $L^p$ class. These functions
may be defined on finite or infinite intervals, or in complex
domains. Another characteristic is the measure of error which can
be the supremum or $L^p$ norm, etc., with a weight usually
introduced in the case of non-compact domains.

But what are the available data to construct such approximations?
From a practical point of view, this is a crucial question. For
example, convolution integrals are very useful in proving
Jackson's theorem, but to actually calculate them one needs
complete knowledge of the function. Obviously, this precludes any
numerical application. In practice, we are generally given a {\it
discrete} set of data (like function values at certain points)
from which we wish to reconstruct the function.

Another feature which is important to us is that we prefer to
construct {\it linear} operators since they are easier to handle.
For example, the so-called ``best approximation" is not at all
``best" for our purposes, since (except in inner product spaces)
it represents a non-linear operator which may be very difficult to
calculate.

A natural candidate that satisfies the above requirements is some
kind of interpolation operator. It is based on a discrete set of
data, and additionally gives a zero error at an increasing number
of well defined points. Depending on the nature of the
interpolation, it is relatively easy to construct such
interpolation operators.

However, there is a significant drawback to these operators. In
the case of Lagrange interpolation, for example, no matter how we
choose the $n$ nodes we get at least an extra $O(\log n)$ factor
compared to the optimal Jackson order of convergence.  And if we
try to avoid this problematic situation by considering
Hermite--Fej\'er interpolation then, although we obtain the
Jackson order of convergence in some cases, this process is
saturated. That is, it will not give the Jackson order of
approximation, for example, for Lip 1 classes, let alone for
classes of functions with higher order of smoothness.

It transpires that the reason for this negative behavior of
interpolating polynomials is the strict restriction on the degree.
In 1963, at an Oberwolfach conference, Paul Butzer raised the
following question: Is it possible to construct an interpolation
process that gives the Jackson order of convergence? Of course,
this was meant in the sense that we should now allow interpolating
polynomials of degree {\it higher} than Lagrange or
Hermite--Fej\'er interpolation, say polynomials of degree $\le an$
for some constant $a$. It was G\'eza Freud [11] who answered this
question affirmatively by constructing polynomials $p_n$ of degree
at most $4n-3$, interpolating at $\sim n/3$ nodes and
approximating to the Jackson order

$$\displaystyle\max_{|x|\le1}|f(x)-p_n(x)|\le c\,\omega(f,1/n),\qquad|x|\le 1,$$
where $\omega$ is the ordinary modulus of continuity. This was
later improved upon by Freud and V\'ertesi [14] who constructed a
sequence of discrete linear polynomial operators (DLPO, in what
follows) $J_n(f,x)$ of degree at most $4n-2$ which interpolate at
the Chebyshev nodes $\cos((2k-1)/(2n))\pi$, $k=1,\dots,n$, and
provide a Timan type (pointwise) estimate

$$\displaystyle|f(x)-J_n(f,x)|\le c[\omega(f,{\sqrt{1-x^2}\over
n})+\omega(f,{1\over n})],\qquad |x|\le
1.\eqno(1)$$ Freud and Sharma [12, 13] further extended this
result to operators based on general Jacobi nodes, and also
succeeded in decreasing the degree of the polynomial to
$n(1+\varepsilon)$, for an arbitrary $\varepsilon>0$. (Of course,
the $c$ in the above then depends upon $\varepsilon$.)

Freud's work [11] initiated a substantial series of papers
exhibiting constructions of a similar character. R.~B.~Saxena [19]
succeeded in constructing a sequence of DLPOs $J_n^\star$ which
realized the even stronger Telyakovskii--Gopengauz estimate

$$\displaystyle|f(x)-J_n^\star(f,x)|\le c\,\omega(f,{\sqrt{1-x^2}\over
n}),\qquad |x|\le 1.$$

Furthermore, for $2\pi$-periodic continuous functions, O.~Kis and
P.~V\'er\-tesi [17] constructed a very simple interpolating
process: Let $\ell_k(x)$ be the fundamental functions of
trigonometric interpolation based on the equidistant nodes
$x_k=2k\pi/(2n+1),\; k=0,\dots,2n$, i.e.,~let $\ell_k(x)$ be that
trigonometric polynomial of order $n$ for which
$\ell_k(x_j)=\delta_{jk}$, $j,k=0,\dots,2n$, and consider the
sequence of operators

$$U_n(f,x)=\sum_{k=0}^nf(x_k)[4\ell_k(x)^3-3\ell_k(x)^4].$$
Then, for every $2\pi$-periodic continuous functions $f$, we have

$$\max_{x\in{\bf R}}|f(x)-U_n(f,x)|\le c\,\omega(f,{1\over n}),$$
and, for this trigonometric polynomial $U_n(f,x)$ of degree at
most $4n$, evidently

$$U_n(f,x_k)=f(x_k),\qquad k=0,\dots,2n,$$
because of the nice identity

$$\sum_{k=0}^n[4\ell_k(x)^3-3\ell_k(x)^4]=1$$
(cf.~Turetskii [27]). With a proper transformation, Kis and
V\'ertesi [17] also constructed an operator that realized the
Timan estimate (cf.~(1)) and interpolates at some nodes. In that
construction, the limit of the ratio of the degree of the
polynomial and the number of points of interpolation was greater
than 1. Later, O.~Kis and J.~Szabados [16] sharpened this to the
following: {\sl Given $0<\varepsilon\le1$, let
$n\ge20/\varepsilon^2$ and $f \in C[-1,1]$. There exists a
sequence of DLPOs $p_n$ of degree at most $n(1+\varepsilon)$ such
that $p_n$ interpolates $f$ in at least $n$ points and}

$$|f(x)-p_n(x)|\le{13\over\varepsilon^2}\,\omega(f,{\pi\sqrt{1-x^2}\over2n}
), \qquad |x|\le 1.$$ This answers a question raised by
Freud and Sharma [12] about the construction of linear operators
of minimal degree compared to the number of points of
interpolation, that at the same time realize the
Telyakovskii--Gopengauz estimate.

In the algebraic case, the above mentioned operators are based on
nodes of orthogonal polynomials which is not the best choice from
a numerical point of view. Usually, data about an unknown function
are given at {\it equidistant} nodes. Since Lagrange interpolation
at equidistant nodes behaves very poorly (the norm of the operator
grows exponentially with the number of nodes), it was not at all
obvious if convergent operators based on equidistant nodes could
be constructed. On p.~581 of [3], R.~DeVore raised this question
while requesting a Jackson order of convergence. One possible
solution is the following: As an intermediate approximation, we
take the rational functions

$$\displaystyle R_n(f,x)={\displaystyle \sum_{k=-n}^n{f(k/n)\over(x-k/n)^4}
\left/\displaystyle \sum_{k=-n}^n
{1\over(x-k/n)^4}\right.}\;, \eqno(2)$$ based on equidistant nodes. These
are a special case of the so-called Shepard operators which
approximate to Jackson order (for details, see later). Then, to
approximate (2) by polynomials, we can use any standard
constructive operator (for example, the one developed by Bojanic
and DeVore [2]). An easy calculation shows that the moduli of
continuity of $f$ and (2) are equivalent, thus solving the problem
(for details, see [24]). Also, if the operator obtained in this
way is of degree $n$, then for some $c$, $0<c<1$, it is possible
to modify it so that the new operator will interpolate at the
equidistant nodes $k/(cn)$, $k=0,\pm1,\dots,\pm cn$ (cf.~[25],
Theorem 3).

With the exception of Lagrange interpolation, all of the above
mentioned operators are {\it saturated}, i.e.,~the order of
convergence cannot be improved upon beyond a certain limit, even
if the function has better and better structural properties. The
question then becomes: Can we construct discrete linear operators
that are not saturated, and approximate ``close" to the order of
best approximation $E_n(f)$? Of course, by Korovkin's theorem,
such operators cannot be positive.  In the case of $2\pi$-periodic
functions, the well-known de la Vall\'ee Poussin means do have the
desired approximation order, but they are not discrete. In [23],
the following family of operators was constructed. Let $j,n$ be
positive and $k,\ell$ nonnegative integers such that
${1\over2}(jn+kn-k+\ell-1)$ is a nonnegative integer. Set

$$s_{jk\ell n}(t)={\sin{jnt\over2}(\sin{nt\over2})^k(\cos{t\over2})^\ell\over
j( n\sin{t\over2})^{k+1}},$$

$$t_\nu={2\pi\nu\over jn},\qquad\nu=0,\pm1,\pm2,\dots,$$
and let

$$S_{jk\ell n}(g,t)=\sum_{\nu=0}^{jn-1}g(t_\nu)s_{jk\ell n}(t),$$
be a discrete linear operator defined on all $2\pi$-periodic
continuous functions $g$. This is a generalization of many of the
classical kernels and operators. $s_{100n}(t)$ ($n$ odd) is the
Dirichlet kernel, $s_{110n}(t)$ is the Fej\'er kernel,
$s_{310n}(t)$ is the de la Vall\'ee Poussin kernel, and
$s_{130n}(t)$ is the Jackson kernel. Correspondingly,
$S_{100n}(g,t)$ ($n$ odd) and $S_{101n}(g,t)$ ($n$ even) are
ordinary interpolation polynomials, $S_{110n}(g,t)$ is the Jackson
polynomial and $S_{310n}(g,t)$ is the discrete version of the de
la Vall\'ee Poussin operator (see [16] and [23]).  These sequences
of operators converge if the so-called Lebesgue constant

$$\Vert L_{jk\ell n}\Vert=\max_{t\in{\bf R}}\sum_{\nu=0}^{jn-1}|s_{jk\ell n}
(t-t_\nu)|$$ remains bounded as $n\to\infty$ (the other parameters
being considered as constants). Namely, the following typical
error estimate holds:

$$\max_{t\in{\bf R}}|S_{jk\ell n}(g,t)-g(t)|\le(1+\Vert L_{jk\ell n}\Vert)
E_n^T(g),\eqno(3)$$ where $E_n^T(\cdot)$ is the error in the best
approximation by trigonometric polynomials. For some special
choices of the parameters $j,k,\ell$, the following relations
hold:

$$\Vert L_{j10n}\Vert={1\over j}\sum_{\nu=1}^j\cot{2\nu-1\over4j}\pi={2\over
\pi}(\log{8\pi\over j}+\gamma)+\alpha_j$$ where
$\gamma=0.5772\cdots$ is the Euler constant and
$0<\alpha_j<{\pi\over72j^2};$

$$\Vert L_{j11n}\Vert\le{2\over\pi}\log j+2.283\qquad{\rm if}\;j\;{\rm
is\;even,}\;n\;{\rm is\;odd};$$

$$\Vert L_{221n}\Vert\le{2\over\sqrt3};$$

$$\Vert L_{332n}\Vert\le{11\over9}.$$

In all of these cases, (3) yields the order of best approximation.
In addition, evidently

$$S_{jk\ell n}(g,t_\nu)=g(t_\nu), \qquad \nu=1,\dots,jn,$$
i.e.,~the ratio of the order of the operator to the number of
interpolation points is

$${jn+kn+\ell-k+1\over2jn},$$
which is fairly close to the optimal ratio 1/2 (which corresponds
to the minimal degree interpolation) provided $j$ is large
compared with $k$.

The obvious transformation $x=\cos t$ converts these results into
estimates for continuous functions on the interval $[-1,1]$. Let

$$P_{jk\ell n}(f,x)=S_{jk\ell n}(f(\cos t),t)$$
and

$$x_\nu=\cos t_\nu, \qquad \nu=0,\dots,jn-1.$$
This new polynomial operator will have the analogous properties

$$P_{jk\ell n}(f,x_\nu)=f(x_\nu), \qquad \nu=0,\dots,jn-1, $$
and

$$\max_{|x|\le1}|f(x)-P_{jk\ell n}(f,x)|\le(1+\Vert L_{jk\ell n}\Vert)E_n(f),$$
where $E_n(\cdot)$ is the error in the best approximation by
algebraic polynomials of degree at most $n$.

Let the number $r_n=m_n/n$ be the ratio of the degree $m_n$ of a
polynomial operator to the number $n$ of nodes where it
interpolates. For Lagrange interpolation we always have
$\lim_{n\to\infty}r_n=1.$ In the case of (convergent)
Hermite--Fej\'er interpolation, $\lim_{n\to\infty}r_n=2$. A
classical result of Bernstein [1] says that if $\varepsilon>0$ is
arbitrary, then there exists a linear polynomial operator with
bounded norm such that

$$\lim_{n\to\infty}r_n=1+\varepsilon.\eqno(4)$$
In general, it is desirable to investigate the relation between
the norm $\Vert L_n\Vert$ of a (not necessarily linear) polynomial
operator of degree $m_n$ and the number of nodes $n$ where it
interpolates. Generalizing the results just mentioned, it can be
shown that for such operators

$$\lim_{n\to\infty}{\Vert L_n\Vert\over\log{n\over m_n-n+2}}>0.$$
This was proved first for Chebyshev nodes [26], and later
generalized by B.~Shekhtman [20] to an arbitrary system of nodes.
We emphasize that this result is true without assuming linearity
of the operator. It holds, for example, for operators of the form

$$L_n(f,x)=p_n(f,x)+\sum_{k=1}^n[f(x_k)-p_n(f,x_k)]q_k(x),$$
where $p_n$ is the best approximating polynomial of degree $m_n$
and the $q_k$ are polynomials of the same degree such that
$q_k(x_j)=\delta_{kj}$, $j,k=1,\dots,n$. Sequences of operators
having property (4) for general systems of nodes $x_k=\cos t_k$,
$k=1,\dots,n$, and approximating to the order
$O(E_{n(1+\varepsilon)})$ were constructed by Erd\H os, Kro\'o and
Szabados [10] who showed that such approximation is possible if
and only if

$$\limsup_{n\to\infty}{N_n(I_n)\over n|I_n|}\le{1\over\pi}\quad{\rm
whenever}\quad\lim_{n\to\infty}n|I_n|=\infty,$$ and

$$\lim_{n\to\infty}n\cdot\min_{1\le k\le n-1}(t_k-t_{k+1})>0,$$
where $N_n(I_n)$ is the number of $t_k$'s in the interval
$I_n\subset[0,\pi]$.

We have, until now, dealt with polynomial operators. From a
numerical point of view, {\it rational functions} are just as
useful as polynomials, and as we shall see, they enjoy certain
advantages over polynomial operators. We have already briefly
mentioned the operator (2). With an exponent 2 instead of 4, and
in the case of bivariate functions, this operator was originally
introduced by D.~Shepard [21]. His interesting paper went
unnoticed for some time, but his operator was rediscovered in 1976
by J.~Bal\'azs (oral communication). It was then proved in [24]
that for the sequence of operators

$$\displaystyle R_n(f,x)={\displaystyle\sum_{k=-n}^n{f(k/n)\over|x-k/n|^p}
\left/\sum_{k=-n}^n
{1\over|x-k/n|^p}\right.}\;,\eqno(5)$$ in the case $p=4$, we have

$$\max_{|x|\le1}|f(x)-R_n(f,x)|\le c\,\omega(f,{1\over n}),$$
i.e.,~a Jackson order of convergence. The major advantage of this
sequence of discrete linear operators is that it is based on
equidistant nodes (and not on the roots of some orthogonal
polynomials). In addition, these operators (for $p=2$) have the
properties

$$R_n(f,k/n)=f(k/n),\quad R_n^\prime(f,k/n)=0, \qquad k=0,\dots,\pm n,$$
 i.e.,~they behave like Hermite--Fej\'er interpolating polynomials.

Bal\'azs' rediscovery of Shepard's operators initiated a series of
papers dealing with further problems and generalizations. We
mention only a few significant papers and results.  G.~Somorjai
[22] proved that, for $p>2$,

$$\max_{|x|\le1}|f(x)-R_n(f,x)|\le{c\over n}\quad \Leftrightarrow\quad f\in{\rm
Lip\,1}, $$ and

$$\max_{|x|\le1}|f(x)-R_n(f,x)|=o\left({1\over n}\right)\quad \Leftrightarrow\quad
f={\rm const.},$$ i.e.,~the saturation property of the operator
$R_n$. The important (and much more difficult) case $p=2$ was
handled by Della Vecchia, Mastroianni and Totik [4].  They proved
that for $p=2$ the trivial class (for which the approximation
order is $o(1/n)$) is the set of constant functions, and the
saturation class (for which the order of approximation is
$O(1/n)$) is contained in $\cap_{\alpha<1}{\rm Lip}\,\alpha$. They
also proved the surprising fact that, for arbitrary {\it
subsequences} of the Shepard operators, arbitrarily fast
convergence may be achieved for some nonconstant functions.

The case $0<p\le2$ was handled by X.~L.~Zhou [28]. With the
notation

$$\displaystyle T_{\varepsilon,p}f:=\max_{0\le x\le1}\left|
\int_{0\le
t\le1,\,|t-x|\ge\varepsilon}{f(t)-f(x)\over|t-x|^p}\,dt\right|,$$
he proved that the saturation order is $n^{1-p}$ or $1/\log n$,
and the saturation class is

$$\displaystyle\{f: f\in\,{\rm Lip}\,(p-1),\;
\sup_{\varepsilon>0}\,T_{\varepsilon,p}f<\infty\}$$ or

$$\displaystyle\{f: \omega(f,t)=O(|\log t|^{-1}),\;
\sup_{\varepsilon>0}\,T_{\varepsilon,1}f<\infty\}$$ according as
$1<p\le2$ or $p=1$, respectively. He also showed that the case
$0<p<1$ is not interesting, since then

$$\displaystyle\max_{0\le x\le1}|f(x)-R_n(f,x)|\to0\qquad\Leftrightarrow
\qquad f={\rm const.}$$

Many of the above mentioned results have been generalized to
weighted approximation (see [6]--[9], [18]), as well as to
general systems of nodes (see [5]).

\References


\noindent\hangindent5mm [1] S.~Bernstein, Sur une modification de
la formula d'interpolation de Lagrange, {\it Zap.~Hark.~Mat. Tov.}
{\bf5} (1932), 49--57.

\noindent\hangindent5mm [2] R.~Bojanic and R.~DeVore, A proof of
Jackson's theorem, {\it Bull.~Amer. Math.~Soc.} {\bf75} (1969),
364--367.

\noindent\hangindent5mm [3] P.~L.~Butzer and B.~Sz.~Nagy (eds.),
{\it Linear Operators and Approximation II,} Proceeding of the
Conference held in Oberwolfach, ISNM Vol.~25, Birkh\"auser (1974),
585 pp.

\noindent\hangindent5mm [4] B.~Della Vecchia, G.~Mastroianni and
V.~Totik, Saturation of the Shepard operators, {\it Approx.~Theory
and its Appl.} {\bf6} (1990), 76--84.

\noindent\hangindent5mm [5] B.~Della Vecchia and G.~Mastroianni,
Pointwise simultaneous approximation by rational operators, {\it
J.~Approx.~Theory} {\bf65} (1991), 140--150.

\noindent\hangindent5mm [6] B.~Della Vecchia, G.~Mastroianni and
J.~Szabados, Weighted uniform approximation on the semiaxis by
rational operators, {\it J.~of Inequalities and Appl.} {\bf4}
(1999), 241--264.

\noindent\hangindent5mm [7]  B.~Della Vecchia, G.~Mastroianni and
J.~Szabados, Approximation with exponential weights in $[-1,1]$,
{\it J.~Math.~Anal.~Appl.} {\bf272} (2002), 1--18.

\noindent\hangindent5mm [8]  B.~Della Vecchia, G.~Mastroianni and
J.~Szabados, Weighted approximation of functions with inner
singularities by exponential weights in $[-1,1]$, {\it
Numer.~Funct.~Anal. and Optimiz.} {\bf 24} (2003), 181--194.

\noindent\hangindent5mm [9]  B.~Della Vecchia, G.~Mastroianni and
J.~Szabados, Weighted approximation of functions with endpoint or
inner singularities by iterated exponential weights, {\it East
J.~Approx.} {\bf9} (2003), 215--227.

\noindent\hangindent5mm [10] P.~Erd\H os, A.~Kro\'o and
J.~Szabados, On convergent interpolatory polynomials, {\it
J.~Approx.~Theory} {\bf58} (1989), 232--241.

\noindent\hangindent5mm [11] G.~Freud, \"Uber ein Jacksonsches
Interpolationsverfahren, in: {\it On Approximation Theory,}
Proc. Conference in Oberwolfach (eds.~P.~L.~Butzer and
J.~Korevaar), ISNM Vol.5 (Birkh\"auser, 1963), pp.~227--232.

\noindent\hangindent5mm [12] G.~Freud and A.~Sharma, Some good
sequences of interpolation polynomials, {\it Canad.~J.~Math.}
{\bf26} (1974), 233--246.

\noindent\hangindent5mm [13] G.~Freud and A.~Sharma, Addendum:
Some good sequences of interpolation polynomials, {\it
Canad.~J.~Math.} {\bf29} (1977), 1163--1166.

\noindent\hangindent5mm [14] G.~Freud and P.~V\'ertesi, A new
proof of A.~F.~Timan's approximation theorem, {\it Studia
Sci.~Math.~Hungar.} {\bf2} (1967), 403--414.

\noindent\hangindent5mm [15] D.~Jackson, A formula of
trigonometric interpolation, {\it Rendiconti del Circolo
Mathematico di Palermo} {\bf 37} (1914), 1--37.

\noindent\hangindent5mm [16] O.~Kis and J.~Szabados, On some de la
Vall\'ee Poussin type discrete linear operators, {\it Acta
Math.~Hungar.} {\bf47} (1986), 239--260.

\noindent\hangindent5mm [17] O.~Kis and P.~V\'ertesi, On a new
process of interpolation (in Russian), {\it
Ann.~Univ.~Sci.~Budapest., Sectio Math.} {\bf10} (1967), 117--128.

\noindent\hangindent5mm [18] G. Mastroianni and J.~Szabados,
Bal\'azs--Shepard operators on infinite intervals. II, {\it
J.~Approx.~Theory} {\bf90} (1997), 1--8.

\noindent\hangindent5mm [19] R.~B.~Saxena, Approximation of
continuous functions by polynomials, {\it Studia
Sci.~Math.~Hungar.} {\bf8} (1973), 437--446.

\noindent\hangindent5mm [20] B.~Shekhtman, On the norms of
interpolating operators, {\it Israel J.~Math.} {\bf64} (1988),
39--48.

\noindent\hangindent5mm [21] D.~Shepard, A two dimensional
interpolation function for irregularly spaced data, {\it
Proc.~23rd National Conference ACM} (1968), 517--523.

\noindent\hangindent5mm [22] G.~Somorjai, On a saturation problem,
{\it Acta Math.~Acad.~Sci.~Hungar.} {\bf32} (1978), 377--381.

\noindent\hangindent5mm [23] J.~Szabados, On an interpolatory
analogon of the de la Vall\'ee Poussin means, {\it Studia
Sci.~Math. Hungar.} {\bf9} (1974), 187--190.

\noindent\hangindent5mm [24] J.~Szabados, On a problem of
R.~DeVore, {\it Acta Math.~Acad.~Sci.~Hungar.} {\bf27} (1976),
219--223.

\noindent\hangindent5mm [25] J.~Szabados, On some convergent
interpolatory polynomials, in: {\it Fourier Analysis and
Approximation Theory,} Coll.~Math.~Soc.~J\'anos Bolyai, Vol.~19
(Budapest, 1976), pp.~805--815.

\noindent\hangindent5mm [26] J.~Szabados, On the norm of certain
interpolating operators, {\it Acta Math.~Hungar.} {\bf55} (1990),
179--183.

\noindent\hangindent5mm [27] A.~H.~Turetskii, On some extremal
problems in the theory of interpolation, in: {\it Investigations
in Contemporary Problems of Constructive Theory of Functions},
Akad.~Nauk USSR (Baku, 1965), pp.~220--232 (in Russian) (Minsk,
1960).

\noindent\hangindent5mm [28] X.~L.~Zhou, The saturation class of
Shepard operators, {\it Acta Math.~Hungar.} {\bf80} (1998),
293--310.

{

\bigskip\obeylines
J.~Szabados
Alfr\'ed R\'enyi Institute of Mathematics
P.~O.~Box 127
Budapest H-1364
Hungary
{\tt szabados@renyi.hu}

}


\bye